\documentclass[11pt]{article}
\usepackage{}
\usepackage[total={6in, 9in}]{geometry}
 \usepackage{amsfonts}
\usepackage{amsthm}
\usepackage{amssymb}
\usepackage{amsmath,amsfonts}
\usepackage{mathrsfs}
\usepackage{cases}
\usepackage{latexsym,bm}
\usepackage{indentfirst}
\usepackage{color}
\usepackage{ifpdf}
\usepackage{graphicx}
\usepackage{psfrag}
\usepackage{dsfont}
\usepackage{enumerate}
\usepackage[font=small,labelfont=bf, width=.85\textwidth]{caption}
\usepackage{pgf,tikz,pgfplots}
\usepackage{tikz-3dplot}

\newtheorem{THM}{\textbf{Theorem}}

\newtheorem{LEM}[THM]{\textbf{Lemma}}
\newtheorem{CLA}{\textbf{Claim}}
\newtheorem{CON}{\textbf{Conjecture}}

\newcommand{\pf}{\noindent\textbf{Proof}.\quad}

\newcommand{\ve}{\varepsilon }

\DeclareMathOperator{\df}{def}
\DeclareMathOperator{\odd}{odd}
\newcommand{\pbar}{\overline{\varphi}}
\newcommand{\sbar}{\overline{\psi}}

\begin{document}
\title{Total coloring graphs with large maximum degree}
\author{
	Aseem Dalal\footnote{Indian Institute of Technology Delhi, Department of Mathematics, Delhi India.
		Email: {\tt aseem.dalal@gmail.com}.}
	\qquad
	Jessica McDonald\footnote{Auburn University, Department of Mathematics and Statistics, Auburn U.S.A.
		Email: {\tt mcdonald@auburn.edu}.
		Supported in part by Simons Foundation Grant \#845698  }
	\qquad
	Songling Shan\footnote{Auburn University, Department of Mathematics and Statistics, Auburn U.S.A.
		Email: {\tt szs0398@auburn.edu}.
		Supported in part by NSF grant DMS-2345869.}
}

\date{}
\maketitle

\begin{abstract} We prove that for any graph $G$, the total chromatic number of $G$ is at most $\Delta(G)+2\lceil \frac{|V(G)|}{\Delta(G)+1} \rceil$.
This saves one color in comparison with a result of Hind from 1992. In particular, our result says that if $\Delta(G)\ge \frac{1}{2}|V(G)|$, then $G$ has a total colouring using at most $\Delta(G)+4$ colors.  When $G$ is regular and has a sufficient number of vertices, we can actually save an additional two colors. Specifically, we prove that for any  $0<\ve <1$, there exists $n_0\in \mathbb{N}$ such that: if $G$ is an $r$-regular graph on $n \ge n_0$ vertices with $r\ge \frac{1}{2}(1+\ve) n$, then $\chi_T(G) \le \Delta(G)+2$. This confirms the Total Coloring Conjecture for such graphs $G$.
\end{abstract}

\section{Introduction}

In this paper, all graphs are simple unless otherwise explicitly stated. For terms not defined here, we follow \cite{WestText}.

A \emph{total coloring} of a graph $G$ is an assignment of colors to the vertices and edges of $G$ such that adjacent vertices receive different colors, adjacent edges receive different colors, and if any edge $e$ is incident to any vertex $v$, then $e, v$ also receive different colors.
The  \emph{total chromatic number} of $G$, denoted $\chi_T(G)$, is the minimum number of colors needed to totally color $G$. Since a vertex of maximum degree must receive different colors from all its incident edges, every graph $G$ satisfies $\chi_T(G) \ge \Delta(G)+1$. On the other hand,  Behzad~\cite[Conjecture 1, p.44]{Behzad-total-coloring},  and independently Vizing~\cite{Vizing-total-coloring},  proposed the following conjecture.

\begin{CON}[Total Coloring Conjecture]\label{con:total-coloring}
Every graph $G$ has  $\chi_T(G) \le \Delta(G)+2$.
\end{CON}

There has been extensive work towards the Total Coloring Conjecture; even the planar case, which is not yet completely verified, has a significant literature. An early survey of total coloring can be found in the 1995 book of Jensen and Toft \cite{JT} (Section 4.9); see also the survey paper of Geethaa,  Narayananb, and Somasundaram ~\cite{MR4679850} for a more modern account.

In this paper we focus on the Total Coloring Conjecture for graphs with large maximum degree. For any graph $G$ with sufficiently large maximum degree, Molloy and Reed~\cite{Reed-Molly} showed that $\chi_T(G) \le \Delta(G)+10^{26}$. When the maximum degree is large in comparison with $|V(G)|$, this can be improved. In particular, Hind~\cite{Hind}  proved that  $\chi_T(G) \le \Delta(G)+2\lceil  \frac{|V(G)|}{\Delta(G)} \rceil +1$ for any graph $G$. Our first result in this paper is to improve Hind's bound by (at least) one.

\begin{THM}\label{main}
Every graph $G$  satisfies  $\chi_T(G) \le \Delta(G)+2\left\lceil  \frac{|V(G)|}{\Delta(G)+1} \right\rceil$.
\end{THM}

When $\Delta(G) \ge \frac{3}{4} |V(G)|$, Hilton and Hind~\cite{MR1226136}  have already confirmed  Conjecture \ref{con:total-coloring}. However, for graphs with $\Delta(G) < \frac{3}{4}|V(G)$, Theorem \ref{main} offers an improved approximation to the Total Coloring Conjecture. In particular, for all graphs with $\Delta(G)\ge \frac{1}{2}|V(G)|$, Theorem \ref{main} says that $G$ has a total coloring using at most $\Delta(G)+4$ colors.

When the maximum degree of $G$ is large in comparison to $|V(G)|$ \emph{and} $G$ is regular, more can be said. All $r$-regular graphs $G$ with  $r > \frac{2}{3}|V(G)|+\frac{23}{6}$  are known to satisfy the Total Coloring Conjecture (Xie and He~\cite{MR2170126} and Xie and Yang~\cite{MR2530186}); this result is built on earlier work by Chetwynd, Hilton and Zaho \cite{MR1136434}  and also Chew~\cite{MR1703308}. Our second result in this paper is to show that for graphs on a sufficient number of vertices, this fraction of $\tfrac{2}{3}$ can be lowered all the way to $\tfrac{1}{2}$.

\begin{THM}\label{thm:main0}
For any  $0<\ve <1$, there exists $n_0\in \mathbb{N}$ such that: if $G$ is an $r$-regular graph on $n \ge n_0$ vertices such that $r\ge \frac{1}{2} (1+\ve) n$, then $\chi_T(G)  \le  \Delta(G)+2$.
\end{THM}

The remainder of this paper is organized as follows.  In Section 2, we use a Vizing-fan type structure to prove two specific \emph{edge-coloring} results: one for a hypergraph that is nearly a graph, and one for a graph that is nearly bipartite. In Section 3 we use the hypergraph edge-colouring result, along with a classic result of Hajnal and Szemeredi~\cite{MR0297607} on \emph{equitable vertex-colorings}, to quickly prove Theorem \ref{main}. The nearly-biparitite edge-coloring result is needed (along with the hypergraph edge-coloring result) for our proof of Theorem \ref{thm:main0}. Section 4 contains some other needed preliminaries for Theorem \ref{thm:main0}, divided into two subsections on matchings (in the complement of a regular graph), and edge-colorings (which are somewhat equitable). At the start of Section 5 we state one previous result by the third author, about a vertex partition that will be essential for us, and then we present our proof of Theorem \ref{thm:main0}.

\section{Two edge-coloring results}

A \emph{hypergraph} $H$ is a set of vertices $V(H)$ and a set of edges $E(H)$, each of which is a subset of $V(H)$ of size at least two. In this paper we assume that all hypergraphs are simple, meaning that any two edges intersect in at most one vertex. A hypergraph $H$ is said to have rank $c$ if the largest set in $E(H)$ has size $c$; if all edges in $H$ have size 2 then $H$ is a graph. Given a hypergraph $H$, a vertex $v$ has \emph{degree}, denoted $d_H(v)$,  equal to the number of edges of $H$ it is a member of; $\Delta(H)$ denotes the maximum degree of $H$. A \emph{matching} $M$ in a hypergraph $H$ is a set of edges, none of which share any common vertices; a vertex is \emph{saturated} by $M$ if it is contained in an edge of $M$; $V(M)$ denotes the set of all vertices in $V(H)$ that are saturated by $M$, and; $M$ is \emph{perfect} if $V(H)=V(M)$.  Given an integer $k\ge 1$, a \emph{$k$-edge-coloring} of a hypergraph $H$ (with at least one edge) is an assignment of the colors $1,  \ldots, k$ to the edges of $H$ such that two edges sharing a common vertex must receive different colors. For a vertex $x$ in a hypergraph $H$, by $\delta_H(x)$ we mean the set of all edges containing $x$ in $H$; by $N_H(x)$ we mean the set of all vertices in an edge with $x$ in $H$.

Given a hypgergraph $H$, a partial $k$-edge-coloring of $H$ is an edge-coloring of some sub-hypergraph $H'\subseteq H$. For $\varphi$ any partial $k$-edge-coloring of $H$, and any vertex $v\in V(H)$, denote by $\pbar(v)$ the set of colors \emph{missing} at $v$, that is, the set of all colors in $\{1, 2, \ldots, k\}$ that are not used on an edge incident to $v$ under $\varphi$. For a vertex set $X\subseteq V(H)$,  define  $\pbar(X)=\bigcup _{v\in X} \pbar(v)$ to be the set of missing colors of $X$. Given any pair of distinct colors $\alpha, \gamma \in\{1, \ldots, k\}$, and given a vertex $v\in V(H)$ missing at least one of $\alpha, \gamma$ under $\varphi$, we can consider the \emph{maximal $(\alpha, \gamma)$-chain  $P$ beginning at $v$}, that is,  the component containing $v$ of
the sub-hypergraph induced by edges colored $\alpha$ or $\gamma$ under $\varphi$. The \emph{end-vertices} of such a chain are the vertices that have degree one in the chain.
If $H$ is a graph, then a maximal $(\alpha, \gamma)$-chain is a maximal $(\alpha, \gamma)$-alternating path, and its end-vertices are the endpoints of the path. Given such a maximal chain  $P$, note that by switching the colors $\alpha, \gamma$ on $P$ we get a new partial $k$-edge-coloring $\varphi'$ of $H$ where exactly the same set of edges in $H$ have been colored; we denote this $\varphi'$ by $\varphi/P$.

Let $G$ be a graph and let $\varphi$ be a partial $k$-edge-coloring of $G$.  A \emph{multifan} centered at $r$ with respect to $\varphi$ is a sequence $F_\varphi(r,s_0: s_p)=(r, e_0, s_0, e_1, s_1, \ldots, e_p, s_p)$ with $p\geq 0$ consisting of  distinct edges $e_0, \ldots, e_p$  with $e_i=rs_i$ for all $i$,  where $e_0$ is left uncolored under $\varphi$,
and \emph{for every edge $e_i$ with $i\in \{1, \ldots, p\}$,  there exists  $j\in \{0, \ldots, i-1\}$ such that
		$\varphi(e_i)\in \pbar(s_j)$.}

The above definition of a multifan  first appeared in the book of Stiebitz et al.~\cite{StiebSTF-Book} as a generalization of a \emph{Vizing-fan}, which Vizing used to prove his eponymous theorem.
Vizing required his fans to be simple, while multifans are typically allowed to have parallel edges, although that difference is not relevant for us here, since all graph (and hypergraphs) are simple in this paper. The other difference in the structures is important for us here however: in the italicized part of the definition above, Vizing required $j=i-1$.
We don't require this, but we will however find it useful to consider multifans which have this property. To this end, first note that if $F_\varphi(r,s_0: s_p)$ is a multifan, then for any integer $p^*\in \{0,  \ldots, p\}$, $F_\varphi(r,s_0: s_{p^*})$ is also a multifan.
A subsequence $(s_0=s_{\ell_0}, s_{\ell_1},s_{\ell_2}, \ldots, s_{\ell_t})$ with the property that $\varphi(e_{\ell_i})= \alpha\in \pbar(s_{\ell_{i-1}})$ for each $i\in \{1,\ldots, t\}$, is called a \emph{linear sequence} of  $F_\varphi(r,s_0: s_{p})$.
Given such a linear sequence, and given $h \in \{1, \ldots, t\}$, we \emph{shift from $s_{\ell_h}$ to $s_0$} by recoloring edge $e_{\ell_{i-1}}$ with $\varphi(e_{\ell_{i}})$ for all $i\in\{1, \ldots, h\}$.
Note that the result is a partial $k$-edge-coloring where $e_0$ is colored and $e_{\ell_h}$ is not. These linear subsequences and shifts have been implicitly used in many papers (see eg. \cite{MR4694336}).

The first of the two edge-coloring results we will prove in this section is as follows.

\begin{LEM}\label{bipartite-matching-extension}
	Let $H$ be a hypergraph with rank $c$ and suppose that there exists a matching $M$ of $H$ such that all edges in $E(H)\setminus M$ have size two. Suppose further that there exists a vertex $x\in V(H)$ with $N_H(x)\subseteq V(H)\setminus V(M)$.
If $|M|+|\delta_H(x)| \le \Delta(H)+2c-1,$ then $H$ has a $(\Delta(H)+2c-1)$-edge-coloring where every edge of $M\cup\delta_H(x)$ receives a different color.
\end{LEM}

Note that the $H$ in Lemma \ref{bipartite-matching-extension} is a hypergraph, but it is nearly a graph (all but the matching edges are of size two). Our approach will be to find a multifan within the edges of size two, and to consider linear subsequences and shifts of this multifan, but also to use maximal $(\alpha, \gamma)$-chains in the hypergraph which may indeed contain edges of larger sizes.

\proof[Proof] \emph{(Lemma \ref{bipartite-matching-extension})}  Let $\Delta=\Delta(H)$ and let $M^*=M\cup \delta_H(x)$. By assumption,
$$|M^*|=|M|+|\delta_H(x)| \le \Delta+2c-1.$$
Let $\varphi$ be a partial  $(\Delta+2c-1)$-edge-coloring of $H$ with the following properties:
\begin{enumerate}
\item[(P1)] every edge in $M^*$ receives a different color, and;
\item[(P2)] subject to (P1), $\varphi$ assigns colors to as many edges as possible.
\end{enumerate}
If every edge in $H$ is colored by $\varphi$ then we are done, so we may assume there is some edge $e_0$ that is uncolored by $\varphi$. Since $e_0$ is uncolored it is not in $M^*$, so by assumption $e_0=uv$ for some $u, v\in V(H)\setminus\{x\}$. Let $F$ be a  multifan centered at $u$, which is maximal subject to only containing edges of size two.

Since $u\neq x$, $u$ is incident with at most one edge from $M^*$. Since $\varphi$ has $\Delta+2c-1$ colors (and $uv$ is uncolored), we also know that $|\pbar(v)| \ge 2c\geq 4$, and all of these  colors must be present at $u$ (otherwise $\varphi$ could be extended to $e_0=uv$).
So there are certainly linear sequences of $F$ containing no edge of $M^*$. We let $F_1\subseteq F$ be a multifan induced by all the linear sequences of $F$ that don't contain edges of $M^*$.
Suppose $F_1=(u, uw_0, w_0,  uw_1, w_1, \ldots,  uw_t, w_t)$,  where $w_0=v$.

\begin{CLA}\label{uwi} Suppose $F_2\subseteq F_1$ corresponds to the union of any number of linear sequences of $F$. Then $\pbar(u)\cap \pbar(w_i) =\emptyset$ for any $w_i\in V(F_2)$.
\end{CLA}

\proof[Proof of claim] We already explained the truth of this for $w_0=v$, and can extend this argument to the other $w_i$ as well. To this end, suppose that $\pbar(u)\cap \pbar(w_i) \ne \emptyset$ for some $i$. Then there exists a linear sequence $(u, w_0, w_{i_1}, \ldots, w_{i_s})$ of $F_2$ for which $w_{i_s}=w_i$. Shift in $(u, w_0, w_{i_1}, \ldots, w_{i_s})$
from $w_{i_s}$ to $w_0$ and then color the edge $uw_i$ by a color in $\pbar(u)\cap \pbar(w_i)$.  The resulting partial $(\Delta+2c-1)$-edge-coloring  $\varphi'$ of $H$ still has every edge in $M^*$  colored differently, since we excluded the edges of $M^*$ from $F_1$ (and $F_2$). Hence $\varphi'$ contradicts (P2).
\qed

\begin{CLA}\label{gamma} For any color $\gamma \in \{1, \ldots, \Delta+2c-1\}$,  there are at most  $2c-1$  vertices from $\{w_0, \ldots, w_t\}$ that are all missing $\gamma$ under $\varphi$.
\end{CLA}
	
\proof[Proof of claim]  Suppose to the contrary that there exists such a $\gamma$, and let  $y_1, y_2,  \ldots,  y_{2c}$ be $2c$  distinct vertices of $\{w_0, \ldots, w_t\}$ all of which are missing $\gamma$ under $\varphi$. Suppose, without loss of generality, that $y_1$ is the first vertex
in the order $w_0, w_1, \ldots, w_t$ for which $\gamma\in \pbar(y_1)$.  In particular, this implies that $F_1(u, w_0:y_1)$ does not contain any edge that is colored by $\gamma$.
Let $\alpha\in \pbar(u)$; we know that $\alpha\ne \gamma$ by Claim \ref{uwi}. Consider the set $\mathcal{P}$ of maximal $(\alpha, \gamma)$-chains staring at $u, y_1,\ldots, y_{2c}$, respectively.
At most one edge of $M^*$ is colored by $\alpha$, and at most one edge of $M^*$ is colored by $\gamma$, and any other edge coloured $\alpha$ or $\gamma$ has size two. So, if a chain in $\mathcal{P}$ contains exactly one edge of $M^*$, then it has at most $c$ end-vertices; if it contains exactly two edges of $M^*$, then it can have at most $2c$ end-vertices. Since there are at least $2c+1$ different end-vertices for the chains in $\mathcal{P}$, this means that there is at least one chain $P\in \mathcal{P}$ not containing  any edges from $M^*$. In particular, all edges in this chain will have size two. We let $\varphi'=\varphi/P$ and consider two different cases according to whether or not $u$ is an endpoint of $P$.

\noindent{\bf Case 1}:  \emph{The vertex $u$ is an endpoint of $P$.}

Suppose first that the other end of $P$ is $y_1$. Then in $\varphi'$, $u$ is missing $\gamma$, so in particular no $\gamma$-edges appear in $F_1$ -- if there were any before, they must have been recolored to $\alpha$. However since these edge all occur after $y_1$ in the order of $F_1$, and $y_1$ is missing $\alpha$ in $\varphi'$, $F_1$ is still a multifan in $G$ under $\varphi'$. However now Claim \ref{uwi} is not satisfied for $\varphi'$, leading to a contradiction in our choice of $\varphi$.

We may now assume that the other end of $P$ is not $y_1$. Here, at least $F_1(u, w_0:y_1)$ is a multifan with respect to $\varphi'$ containing no edges of $M^*$. But now the color $\gamma$ is missing at both $u$ and $y_1$ under $\varphi'$, so Claim \ref{uwi} is not satisfied for $\varphi'$, leading to a contradiction in our choice of $\varphi$.

\noindent{\bf Case 2}:  \emph{The vertex $u$ is not an endpoint of $P$.}

In this case we know that $P$ does not contain any edges of $F_1$.

First suppose that $y_1$ is not an end of $P$. Then $F_1$ is still a multifan with respect to $\varphi'$ containing no edges of $M^*$. However then color $\alpha$ is missing at $u$ and at least one other vertex from $\{y_2, y_3, y_4\}$  under $\varphi'$. This means that Claim \ref{uwi} is not satisfied for $\varphi'$, leading to a contradiction in our choice of $\varphi$.

We may now assume that $y_1$ is an end of $P$. Here, at least $F_1(u, w_0:y_1)$ is a multifan with respect to $\varphi'$ containing no edges of $M^*$. But now the color  $\alpha$ is missing at both $u$ and $y_1$ under $\varphi'$, so Claim \ref{uwi} is not satisfied for $\varphi'$, leading to a contradiction in our choice of $\varphi$.
\qed

Claims \ref{uwi} and \ref{gamma} tell us that
\begin{equation}\label{VF1a}
	|\pbar(V(F_1))|  \ge  |\pbar(u)|+ \tfrac{1}{2c-1}\sum_{i=0}^t|\pbar(w_i)|.
\end{equation}
If we let $H'$ be the hypergraph induced by all the colored edges in $H$, then $|\pbar(z)|\geq \Delta+2c-1-d_{H'}(z)$ for all $z\in V(H)$. Since $d_{H'}(z)\leq \Delta$, we know that $|\pbar(z)|\geq 2c-1$ for any $z$. Since the edge $uw_0$ is uncolored, we further know that $|\pbar(w_0)|\geq 2c$. So from (\ref{VF1a}), we get
\begin{equation}\label{pbarLbound}
|\pbar(V(F_1))|  \ge  (\Delta+2c-1-d_{H'}(u))+ \tfrac{1}{2c-1}(2c)+ \tfrac{1}{2c-1}((2c-1)t)>\Delta+2c+t-d_{H'}(u).
\end{equation}
On the other hand, we know that there are exactly $t$ colored edges in $F_1$,  so there are exactly $d_{H'}(u)-t$ coloured edges incident to $u$ that are not included in $F_1$. At most one of these excluded colored edges could be in $M^*$ (since $u\neq x$), while all the rest must have size two (by assumption), and hence under $\varphi$ they must have a color not in $\pbar(\{w_0, \ldots, w_t\})$, by (P2) and by definition of $F_1$. Of course, by virtue of being incident to $u$, none of these excluded colors are in $\pbar(u)$ either. Hence we have found a set of at least $d_{H'}(u)-t-1$ colors from $\{1, \ldots, \Delta+2c-1\}$ that are not in $\pbar(V(F_1)$. So
$$|\pbar(V(F_1))| \leq (\Delta+2c-1)-(d_{H'}(u)-t-1)=\Delta+2c+t-d_{H'}(u),$$
contradicting (\ref{pbarLbound}).
\qed

The second edge-coloring result that we prove in this section is as follows. Note that while Lemma \ref{bipartite-matching-extension} is about a hypergraph that is almost a graph, the following result is about a graph that is nearly a bipartite graph. Despite the different settings, the arguments are very similar.

\begin{LEM}\label{lem:bipartite-matching-extension}
		Let $G_0$ be a  bipartite graph with bipartition $(A, B)$ and let $M$ be a matching in $G_0$. Let $G$ be a graph obtained from $G_0$ by adding a new vertex $x$ and adding some number of edges between $x$
	and $V(G)\setminus V(M)$.  Let $k=\max\{\Delta(G), k_B+1\}$, where $k_B=\max\{d_{G}(v): v\in B\}$.
	 If $k \ge |M|+|\delta_G(x)|$,
	 then $G$ has a $k$-edge-coloring where every edge of $M\cup \delta_G(x)$ receives a different color.
\end{LEM}

\proof[Proof]  Let $\Delta=\Delta(G)$ and let $M^*=M\cup \delta_G(x)$. By assumption, $|M^*| \leq k.$ Let $\varphi$ be a partial  $k$-edge-colouring of $G$ with the following properties:
\begin{enumerate}
	\item[(P1)] every edge in $M^*$ receives a different color, and;
	\item[(P2)] subject to (P1), $\varphi$ assigns colors to as many edges as possible.
\end{enumerate}
If every edge in $G$ is colored by $\varphi$ then we are done, so we may assume there is some edge $uv$ that is uncolored by $\varphi$.  Note that $u, v \ne x$, since otherwise $uv$ would be in $M^*$ and hence be colored under $\varphi$.  Suppose, without loss of generality, that $u\in A$, and let $F$ be a maximal multifan centered at $u$.

We know that $v\in B$ (since $v\neq x$) so $k\geq d_{G}(v)+1$. We also know that $uv$ is uncolored, therefore $|\pbar(v)|\geq 2$. All of these colors must be present at $u$ (otherwise $\varphi$ could be extended to $e_0=uv$). On the other hand, since $u\neq x$, $u$ is incident with at most one edge from $M^*$. So there are certainly linear sequences of $F$ containing no edge of $M^*$. We let $F_1\subseteq F$ be a multifan induced by all the linear sequences of $F$ that don't contain edges of $M^*$.  Suppose $F_1=(u, uw_0, w_0,  uw_1, w_1, \ldots,  uw_t, w_t)$,  where $w_0=v$. In particular, note that $x\neq w_i$ for any $i$, since its incident edges are all in $M^*$.

\begin{CLA}\label{uwi2} Suppose $F_2\subseteq F_1$ corresponds to the union of any number of linear sequences of $F$. Then $\pbar(u)\cap \pbar(w_i) =\emptyset$ for any $w_i\in V(F_2)$.
\end{CLA}

\proof[Proof of claim] This argument is identical to that of Claim \ref{uwi} within the proof of Lemma \ref{bipartite-matching-extension}.
\qed

As $G_0=G-x$ is bipartite,  for any  distinct $\alpha, \beta  \in \{1,\ldots, k\}$ with $\alpha \in \pbar(u), \beta \in \pbar(w_i)$
for some $i\in \{0,\ldots, t\}$,  if there is an  $(\alpha,\beta)$-alternating path joining $u, w_i$, then it creates an odd cycle with $uw_i$, and hence we know  that the path must contain $x$ (as an internal vertex). Moreover, this means that the path contains edges from $M^*$.
Similarly, if $\beta \in \pbar(w_i)\cap \pbar(w_j)$ for
distinct $i,j\in \{1,\ldots, t\}$, and there is an $(\alpha,\beta)$-alternating path joining $w_i$ and $w_j$, then this creates an odd cycle with $uw_i, uw_j$, and hence we know that the path must contain $x$ (as an internal vertex). Again, we also get that such a path contains edges from $M^*$.

\begin{CLA}\label{gamma2} For any   $\gamma \in \{1,\ldots, k\}$,  there is  at most  one vertex  from $\{w_0, \ldots, w_t\}$ that is  missing $\gamma$ under $\varphi$.
\end{CLA}

\proof[Proof of claim]  Suppose to the contrary that $y_1, y_2$ are two  distinct vertices of $\{w_0, \ldots, w_t\}$, both of which are missing $\gamma$ under $\varphi$. Suppose, without loss of generality, that $y_1$ is the first vertex
in the order $w_0, w_1, \ldots, w_t$ for which $\gamma\in \pbar(y_1)$. In particular, this implies that $F_1(u, w_0:y_1)$ does not contain any edge that is colored by $\gamma$.

Let $\alpha\in \pbar(u)$; we know that $\alpha\ne \gamma$ by Claim \ref{uwi2}. Consider the set of  maximal $(\alpha, \gamma)$-alternating paths staring at $u, y_1, y_2$.  If these three paths are all disjoint,
then one of them, say $P$,  does not contain any edge of $M^*$ (since there are at most two edges of $M^*$ with the colors $\alpha, \gamma$). On the other hand, suppose that two of the three paths are in fact the same path $P^*$. By our discussion prior to this claim, we know that $P^*$ must contain $x$ as an internal vertex, and hence there are two $M^*$-edges colored $\alpha, \gamma$ on the path $P^*$. But then again, among the three paths, we find one $P$ that does not contain any edges of $M^*$.

We let $\varphi'=\varphi/P$ and consider two different cases according to whether or not $u$ is an endpoint of $P$. Note that in either case, we have ensured that $P$ doesn't contain any edges of $M^*$, does not contain $x$, and only one of its endpoints is from $F_1$.

\noindent{\bf Case 1}:  \emph{The vertex $u$ is an endpoint of $P$.}

Since the other endpoint of $P$ is outside $F_1$, $F_1(u, w_0:y_1)$ is a multifan with respect to $\varphi'$ containing no edges of $M^*$. But now the color $\gamma$ is missing at both $u$ and $y_1$ under $\varphi'$, so Claim \ref{uwi2} is not satisfied for $\varphi'$, leading to a contradiction in our choice of $\varphi$.

\noindent{\bf Case 2}:  \emph{The vertex $u$ is not an endpoint of $P$.}
In this case we know that $P$ does not contain any edges of $F_1$, and our proof is identical to case 2 of Claim \ref{gamma} in Lemma \ref{bipartite-matching-extension}.
\qed

Claims \ref{uwi2} and \ref{gamma2} tell us that
\begin{equation}\label{VF1}
	|\pbar(V(F_1))|  \ge  |\pbar(u)|+ \sum_{i=0}^t|\pbar(w_i)|.
\end{equation}
If we let $G'$ be the graph induced by all the colored edges in $G$, then $\pbar(z)\geq k-d_{G'}(z)$ for all $z\in V(G)$.
We further argued that $|\pbar(w_0)|\geq 2$ earlier in this proof, and the same argument (less being incident to an uncolored edge) tells us that $|\pbar(w_i)|\geq 1$ for all other $i$. So from (\ref{VF1}), we get
\begin{equation}\label{pbarLbound1}
	|\pbar(V(F_1))|  \ge k-d_{G'}(u)+ 2+ t.
\end{equation}

On the other hand, we know that there are exactly $t$ colored edges in $F_1$,  so there are exactly $d_{G'}(u)-t$ colored edges incident to $u$ that are not included in $F_1$. At most one of these excluded colored edges could be in $M^*$ (since $u\neq x$), while all the rest must have colors that are not in $\pbar(\{w_0, \ldots, w_t\})$, by (P2) and by definition of $F_1$. Of course, by virtue of being incident to $u$, none of these excluded colors are in $\pbar(u)$ either. Hence we have found a set of at least $d_{G'}(u)-t-1$ colors from $\{1,\ldots,k\}$ that are not in $\pbar(V(F_1))$. So we get
$$|\pbar(V(F_1))| \leq k-(d_{G'}(u)-t-1)< k-d_{G'}(u)+ 2+ t,$$
contradicting~\eqref{pbarLbound1}.
\qed

\section{Proof of Theorem \ref{main}}

A \emph{$k$-coloring} of a graph $G$ is an assignment of the colors $1, 2, \ldots, k$ to the vertices of $G$ such that adjacent vertices receive different colors. Given a $k$-coloring $\varphi$, and given any $i\in\{1, \ldots, k\}$, we call $V_i=\{v\in V(G):\varphi(v)=i\}$ the $i$th color class of $\varphi$; if $|V_i|, |V_j|$ differ by at most one for any $i, j\in\{1,  \ldots, k\}$, then we say that $\varphi$ is \emph{equitable}.

\begin{THM}[Hajnal and Szemer\'edi~\cite{MR0297607}]\label{them:Hajnal–Szemeredi}
	Let $G$ be a graph and let $k\ge \Delta(G)+1$ be any integer. Then $G$ has an equitable $k$-coloring.
\end{THM}

We can now show that Theorem \ref{main} follows quickly from Lemma \ref{bipartite-matching-extension} (and Theorem \ref{them:Hajnal–Szemeredi}).

\proof[Proof of Theorem~\ref{main}]  Let $G$ be a graph with $\Delta=\Delta(G)$.
We will show that $G$ has a total coloring using at most $\Delta+2\lceil  \frac{|V(G)|}{\Delta+1} \rceil$ colors.

By Theorem~\ref{them:Hajnal–Szemeredi}, $G$ has  an equitable ($\Delta$+1)-coloring, say $\varphi_0$.   Let $c=\lceil  \frac{|V(G)|}{\Delta+1} \rceil$.
By assumption, each color class of  $\varphi_0$  has size  $c$ or $c-1$.

We construct a hypergraph $H$ as follows. We let $V(G)= V(H)\cup\{x\}$ for one new vertex $x$. We start with $E(G)\subseteq E(H)$, and for any color class $V_i$ of $\varphi_0$, if $|V_i| \ge 2$, we add $V_i$ as an edge to $H$. If there is at least one color class of $\varphi_0$ that has size one, then for each such color class $V_j$, we add $V_j\cup \{x\}$ as an edge of size two to $H$.

Note that $d_{H}(w)= d_G(w)+1$ for all $w\in V(G)$, and
$d_{H}(x) \le \Delta+1$ as $\varphi_0$ has exactly $ \Delta+1$ distinct color classes. Hence $\Delta(H) = \Delta+1$. We now apply Lemma \ref{bipartite-matching-extension} to $H$ with matching $M=\{V_i: \text{$V_i$ is a color class of $\varphi_0$ with size at least 2}\}$ and vertex $x$. To see that this application is possible first note that, $x$ has been defined as an isolate when $V(H)=V(M)$, and otherwise is in an edge of size two with every vertex in $V(H)\setminus V(M)$. The quantity $2c-1\geq 1$, since $c\geq 1$ (as $|V(G)|\geq \Delta+1$). Finally, the quantity $|M|+|\delta_{H}(x)|$ is exactly the number of color classes of $\varphi_0$, which is exactly $\Delta+1=\Delta(H)$. So from Lemma \ref{bipartite-matching-extension} we get a ($\Delta(H)+2c-1(=\Delta+2c)$)-edge-coloring  $\varphi^*$ of $H$ such that all the edges in $M\cup\delta_{H}(x)$ receive different colors.

Define $\varphi:V(G)\cup E(G) \rightarrow \{1,  \ldots, \Delta+2c\}$ by: $\varphi(e)=\varphi^*(e)$ for all $e\in E(G)$; $\varphi(u)=\varphi^*(V_i)$ for all $V_i\in M$ and $u\in V_i$, and;
  $\varphi(w)=\varphi^*(xw)$ for all  $w\in N_{H}(x)$.  Note that this last assignment ensures that every vertex in $G$ receives a color under $\varphi$, since any vertex not in $V(M)$ is adjacent to $x$ in $H$. In $\varphi$, adjacent edges receive different colors since this is true in $\varphi^*$. Any vertex $w\in V(G)$ gets a different color from all its incident edges in $\varphi$ since in $\varphi^*$, that color was used on an additional edge incident to $w$. Finally, because every edge in $M\cup \delta_{H}(x)$ receives a different color under $\varphi^*$, the only  set of vertices receiving the same colors under $\varphi$ are the vertices from the same color class  (with size at least two) of $\varphi_0$, which is an independent set in $G$.
\qed

\section{Additional Preliminaries}

\subsection{Matchings in the complement of a regular graph}

 The following result can be easily proved by applying Hall's  Marriage Theorem,
 and one can also find a short proof in~\cite[Lemma 9]{MR2993074}.
\begin{LEM}\label{lem:matching-in-bipartite}
	Let $G[X,Y]$ be bipartite graph with $|X|=|Y|=n$.  If $\delta(G) \ge n/2$, then $G$ has a perfect matching.
\end{LEM}

Give a graph $G$ and a set $X\subseteq V(G)$, the \emph{deficiency of $X$} is defined to be the quantity $\df(G)=\odd(G-X)-|X|$, where for any graph $H$, $\odd(H)$ is the number of odd components or components of odd order of $H$.
Note that any matching of $G$ must contain at least $\df(X)$ unsaturated vertices, for any $X\subseteq V(G)$.
The deficiency of $G$ is $\df(G)=\max\{\df(X): X\subseteq V(G)\}$. The Tutte-Berge Formula~\cite{MR100850,MR23048} asserts that the maximum size of a matching in $G$ is equal to $\tfrac{1}{2}(|V(G)|-\df(G))$. For any $X\subseteq V(G)$, define an auxiliary bipartite multigraph $H(X)$ by contracting each component in $G-X$ to a single vertex and deleting the edges of $G[X]$. As part of his new proof of the Tutte-Berge Formula and the Gallai-Edmonds Structure Theorem~\cite{Ed,Gal1,Gal2}, West~\cite{MR2788782} defined $H(X)$ and proved (1)-(2) of the following. Here, by a graph $G$ being \emph{factor-critical}, we mean that  $G$ does not have a perfect matching, but $G-v$ does, for every $v\in V(G)$.

\begin{LEM}
\label{lem:matching-structure}
Let $G$ be a graph and let $X\subseteq V(G)$ be a maximal set with $\df(X)=\df(G)$. Then
	\begin{enumerate}[(1)]
		\item  \emph{(West~\cite[Lemma 2]{MR2788782})}  Every component of $G-X$ is factor-critical.
		\item \emph{(West~\cite[Lemma 3]{MR2788782})} The bipartite multigraph $H(X)$ has a matching saturating  $X$.
		\item  We may choose the matching $M$ in (2) so that if $d$ is the maximum degree in $H$ over all vertices in $X$, then $M$ also saturates all vertices of $V(H(X))\setminus X$ having degree at least $d$ in $H$.
	\end{enumerate}
\end{LEM}

To see Lemma \ref{lem:matching-structure}(3), note first that any bipartite multigraph $H$ with bipartition $(X, Y)$  where $d$ is the maximum degree in $H$ over all vertices in $X$, has a matching saturating all vertices of degree at least $d$ in $Y$  (simply apply Hall's Theorem to $(X, Y')$, where $Y'$ are those vertices of $Y$ having degree at least $d$). Then note that given any
two matchings $M_1, M_2$ of $H$ where $M_1$ saturates  $X_1\subseteq X$ and $M_2$ saturates $Y_2\subseteq Y$, there is a matching $M$ of $H$ that saturates both $X_1$ and $Y_2$ (look at $M_1\cup M_2$ and delete an appropriate set of edges).

We can now prove the following.

\begin{LEM}\label{lem:matching-in-completment}
	Let $G$ be an  $n$-vertex  $r$-regular  graph  for some integer $r\ge 0$. Then $ \overline{G}$ has a matching covering
	at least   $ n-\frac{n}{n-r}$ vertices of $ \overline{G}$.
\end{LEM}

\pf
Note that  $\overline{G}$ is $(n-1-r)$-regular. If  $\overline{G}$ has a perfect matching, then we are done, so suppose not. Let $S\subseteq V(G)$ be a maximal set with $\df(S)=\df(\overline{G})$, as in Lemma \ref{lem:matching-structure}. Then by this lemma,  every component of $\overline{G}-S$ is factor-critical.

Let $D$ of be an odd component of $\overline{G}-S$,  and suppose that $|V(D)| \le n-r-1$. Then the number of edges between $D$ and $S$ in $\overline{G}$ is at least
$$|V(D)|(n-r-1-(|V(D)|-1)).$$
The above function is concave down in $|V(D)|$, and given $1\leq |V(D)| \leq n-r-1$, we  get that the function is always at least $n-r-1$. Thus by Lemma~\ref{lem:matching-structure}(3), $\overline{G}$ has a matching saturating  $S$ and all vertices from the small components of $\overline{G}-S$ (i.e. those of size at most $n-r-1$).
As there are at most  $\frac{n}{n-r}$   components of $\overline{G}-S$ with $|V(D)|  \ge n-r$, and  because each of such component is factor-critical, it
follows that $\overline{G}$ has a matching that covers at least $n-\frac{n}{n-r}$ vertices of $\overline{G}$.
\qed

 \subsection{Equitable edge-colorings}

We used equitable vertex-colorings in the proof of Theorem \ref{main}, and we shall use an analog for edge-coloring for Theorem \ref{thm:main0}. The start of this is the following.

\begin{THM}[McDiarmid \cite{MR300623}]\label{lem:equa-edge-coloring}
 	Let $G$ be a graph with chromatic index $\chi'(G)$. Then for all $k\ge \chi'(G)$, there is a $k$-edge-coloring of $G$  where every color class contains either exactly $\lfloor |E(G)|/k \rfloor$ or exactly $\lceil |E(G)|/k \rceil$ edges.
\end{THM}

Given an edge coloring $\varphi$ of $G$ and a given a color $i$, we denote by $\pbar^{-1}(i)$ the set of all vertices in $G$ which are missing color $i$ under $\varphi$. The following parity lemma will be helpful for us (see eg. Lemma 2.1 in Gr\"{u}newald and Steffen \cite{MR2028248}).

\begin{LEM}\label{lem:parity}
	Let $G$ be a graph and $\varphi$ be a $k$-edge-coloring of $G$ for some integer $k\ge \Delta(G)$.
	Then
	$|\pbar^{-1}(i)| \equiv |V(G)| \pmod{2}$ for every color $i\in \{1,\ldots, k\}$.
\end{LEM}

We now prove the following result about a sort of weakly-equitable of edge-coloring.

\begin{LEM}\label{lem:equitable-coloring-precolored-edges}
	Let $G$ be   a graph  and $F\subseteq E(G)$.
	Suppose that for some integer $k\ge \Delta(G)$,
	$G$ has a $k$-edge-coloring such that all the edges in $F$ receive distinct colors. Then $G$ has
	a $k$-edge-coloring $\varphi$ such that for any  distinct $i, j\in \{1,\ldots, k\}$ it holds that $$\left ||\pbar^{-1}(i)| -|\pbar^{-1}(j)|\right | \le 5,$$ and all edges in $F$ receive distinct colors.
\end{LEM}

\pf  We call a $k$-edge-coloring of $G$ \emph{good} if all the edges in $F$ are assigned distinct colors. Among all good $k$-edge-colorings of $G$, we choose $\varphi$
so that
$$ g_\varphi:= \max\{\left ||\pbar^{-1}(i)| -|\pbar^{-1}(j)|\right | : 1\leq i,j \leq k\}
$$
is minimum.
If $g_\varphi  \le 5$, then we are done, so suppose not. This means there are colors $\alpha, \beta\in \{1,\ldots, k\}$
such that
\begin{equation}\label{gvarphi}
\left ||\pbar^{-1}(\alpha)| -|\pbar^{-1}(\beta )|\right | =g_\varphi \ge 6.
\end{equation}
We  may further assume that $\varphi$ has been chosen so as to minimize the number of pairs $\alpha, \beta$ satisfying (\ref{gvarphi}); let $h_{\varphi}$ be the number of these pairs.

Without loss of generality, suppose that $|\pbar^{-1}(\alpha)| -|\pbar^{-1}(\beta )| \ge 6$. There are at most two edges from $F$ that are colored from $\{\alpha, \beta\}$. Thus there are at most four vertices
from $\pbar^{-1}(\alpha)$ that are endvertices of two $(\alpha, \beta)$-alternating paths
that each contain an edge of $F$. At most $|\pbar^{-1}(\beta)|$
vertices from $\pbar^{-1}(\alpha)$ are endvertices of $|\pbar^{-1}(\beta )|$ $(\alpha, \beta)$-alternating paths involving a vertex from
$\pbar^{-1}(\beta)$. Thus there exist distinct vertices $x,y\in \pbar^{-1}(\alpha)$  occurring in one common $(\alpha, \beta)$-alternating path
that contains no edge of $F$. Let $P$ be the $(\alpha, \beta)$-alternating path with endpoints $x, y$, and let $\psi=\varphi/P$.  As $P$ does not contain any edges of $F$, the edges of $F$ are all still colored  different colors under $\psi$, and hence $\psi$ is a good $k$-edge-coloring of $G$. We claim that $g_\psi \le g_\varphi$ and $h_\psi<h_\varphi$, which yields a contradiction to our choice of $\varphi$.

By our choice of $P$, we know that $|\sbar^{-1}(\alpha)|=|\pbar^{-1}(\alpha)|-2$ and $|\sbar^{-1}(\beta)|=|\pbar^{-1}(\beta)|+2$. Given (\ref{gvarphi}) and the definition of $g_{\varphi}$, for any other color $\gamma \in \{1,\ldots, k\}\setminus \{\alpha,\beta\}$,
we have
$$|\pbar^{-1}(\alpha)|  \ge |\pbar^{-1}(\gamma)|  \ge |\pbar^{-1}(\beta)|.$$  But now, for any pair of colors $\{\alpha',\beta'\}$ satisfying (\ref{gvarphi}), we have
\begin{numcases}{||\sbar^{-1}(\alpha')| -|\sbar^{-1}(\beta' )||=}
	g_\varphi & \text{if $\{\alpha', \beta'\}\cap \{\alpha,\beta\}=\emptyset$},  \nonumber \\
	g_\varphi-2 &  \text{if $\alpha'=\alpha$ and $\beta'\ne \beta$},  \nonumber \\
	g_\varphi-2  & \text{if $\alpha' \ne \alpha$ and $\beta'= \beta$},  \nonumber \\
	g_\varphi-4  & \text{if $\alpha' =\alpha$ and $\beta'= \beta$}.   \nonumber
\end{numcases}
Hence we have $g_\psi \le g_\varphi$ and $h_\psi<h_\varphi$, a contradiction to our choice of $\varphi$.
\qed

\section{Proof of Theorem \ref{thm:main0}}

We will need the following result by the third author.

\begin{LEM}[Shan~{\cite[Lemma 3.2]{SHAN2022429}}]\label{lem:partition}
	There exists a positive integer $n_0$ such that for all $n\ge n_0$ the
	following holds. Let $G$ be a graph on $2n$ vertices, and $N=\{x_1,y_1,\ldots, x_t,y_t\}\subseteq V(G)$, where $ 1\le t\le n$ is an integer.
	Then $V(G)$ can be partitioned into two  parts
	$A$ and $B$ satisfying the properties below:
	\begin{enumerate}[(i)]
		\item  $|A|=|B|$;
		\item $|A\cap \{x_i,y_i\}|=1$ for each $i\in \{1,\ldots, t\}$;
		\item $| d_A(v)-d_B(v)| \le n^{2/3}$ for each $v\in V(G)$.
	\end{enumerate}
	Furthermore, one such partition can be constructed in $O(2n^3 \log_2 (2n^3))$-time.
\end{LEM}

We are now ready to present our proof of Theorem \ref{thm:main0}. We restate the theorem here for convenience.

\setcounter{THM}{1}

\begin{THM}
For any  $0<\ve <1$, there exists $n_0\in \mathbb{N}$ such that: if $G$ is an $r$-regular graph on $n \ge n_0$ vertices such that $r\ge \frac{1}{2} (1+\ve) n$, then $\chi_T(G)  \le  \Delta(G)+2$.
\end{THM}

\proof[Proof of Theorem~\ref{thm:main0}]
Let $n_0$ be chosen to be at least  the lower bound of $n$ in Lemma~\ref{lem:partition} and such that $1/n_0 \ll \ve$.
Let $G$ be an $r$-regular graph on $n\ge n_0$ vertices with $r\ge \frac{1}{2}(1+\ve)n$.
By the result of Hilton and Hind~\cite{MR1226136}  that the Total Coloring Conjecture holds for graphs $G'$ with
$\Delta(G')\ge \frac{3}{4} |V(G')|$, we may assume that $r< \frac{3}{4} n$.

Lemma~\ref{lem:matching-in-completment} tells us that $\overline{G}$  has a matching  $M$  saturating at least $n-\frac{n}{n-r}$
vertices.  As  $r< \frac{3}{4} n$,  we get that $M$ leaves less than $4$ vertices unsaturated. When $n$ is even, this means we can choose $M$ so that it saturates exactly $n-2$ vertices, and when $n$ is odd, we can choose $M$ so that it saturates exactly $n-3$ vertices. When $n$ is odd, the degree-sum formula additionally tells us that $r$ is even. In this case, we may choose a matching $M^*\subseteq M$  saturating exactly $r-2$ vertices of $G$.

\begin{figure}[htb]
    \centering
    \includegraphics[width=15cm]{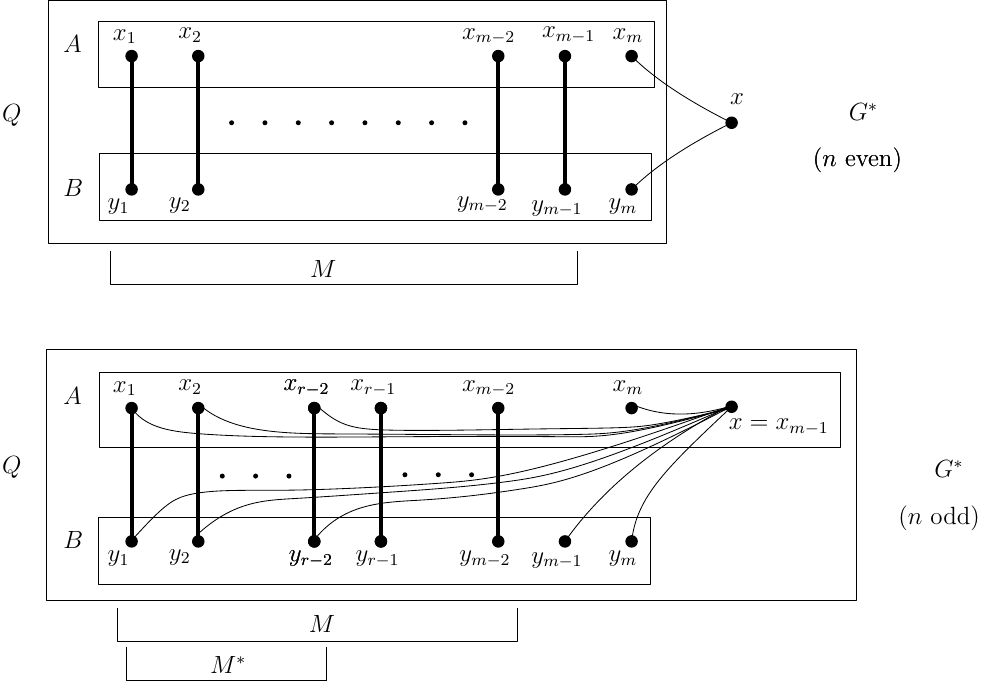}
    \caption{The graph $G^*$ in the cases where $n$ is even and where $n$ is odd. Note the edges of $M$ is bold in each graph, with $|M|=m-1$ in the even case and $|M|=m-2$ in the odd case. Note also that in the odd case we have $M^*\subseteq M$ of size $r-2$, whose endpoints are all joined to $x$.}
    \label{Qpic}
  \end{figure}

We define $G^*$ to be a graph obtained from $G$ by adding the edges of $M$ to $G$, and by adding one new vertex $x$. If $n$ is even, we make $x$ a degree two vertex whose neighbours are the two vertices of $V(G)\setminus V(M)$. When $n$ is odd we make $x$ have degree $r+1$ by joining it to the three vertices of $V(G)\setminus V(M)$ as well as the $r-2$ vertices of $V(M^*)$. See Figure~\ref{Qpic} for a depiction of $G^*$.

Our general aim will be to define an $(r+2)$-edge-coloring of $G^*$ where every edge of $M\cup \delta_{G^*}(x)$ receives a different color. If we can do this (and we will indeed succeed when $n$ is even), then we get a total coloring of $G$ with $r+2$ colors by giving each vertex $w\in V(G)$ the color of the edge from $M\cup \delta_{G}(x)$ that is incident to $w$ in $G^*$. This works, with these auxiliary edges acting as ``placeholders'' for the eventual vertex-coloring, because the only vertices that get the same colour are matched by $M$ (so they are non-adjacent in $G$).

Note that when $n$ is even, $|M \cup \delta_{G^*}(x)|=\tfrac{n-2}{2}+2=\tfrac{n}{2}+1$ which is less than $r+2$ by assumption. On the other hand, when $n$ is odd, $|M \cup \delta_{G^*}(x)|=\tfrac{n-3}{2}+r-1$ which we expect to be strictly larger than $r+2$. So, in the case when $n$ is odd, our stated aim is hopeless, and we will need to adjust. These extra edges from $x$ will however be important at the start of our edge-colouring process (detailed below in five steps). Before the end of our process though, we will be able to safely delete enough of the these edges so as to get our desired result.

Let $m\in\mathbb{Z}^+$ be such that $n=2m$ when $n$ is even and $n=2m-1$ when $n$ is odd. We define the following:
\begin{eqnarray*}
Q&=&G^*-x  \quad \text{if $n$ is even}, \\
Q&=&G^*  \quad \text{if $n$ is odd}, \\
M&=&\{x_1y_1, \ldots, x_{m-1}y_{m-1}\}\quad \text{if $n$ is even},  \\
M&=&\{x_1y_1, \ldots, x_{m-2}y_{m-2}\}\quad \text{if $n$ is odd},  \\
V(Q)\setminus V(M)&=&\{x_{m}, y_{m}\} \quad \text{if $n$ is even},  \\
V(Q)\setminus V(M)&=&\{x_{m-1}, y_{m-1}, x_{m}, y_{m}\} \quad \text{if $n$ is odd, where $x_{m-1}=x$}.
\end{eqnarray*}
Applying Lemma~\ref{lem:partition} on $Q$ and $N=\{x_1,y_1, \ldots, x_{m}, y_{m}\}$,
we obtain a partition $(A, B)$ of $V(Q)$ satisfying the following properties:
\begin{enumerate}[(P1)]
	\item  $|A|=|B|$;
	\item $|A\cap \{x_i,y_i\}|=1$ for each $i\in \{1,\ldots, m\}$;
	\item $\big| d_A(v)-d_B(v)\big| \le m^{2/3}$  for each $v\in V(Q)$.
\end{enumerate}
By renaming vertices if necessary, we assume that $x_i\in A$ and $y_i\in B$ for each $i\in\{1,\ldots, m\}$, with $M^*=\{x_1y_1, \ldots, x_{r-2}y_{r-2}\}$.
Let $Q_A=Q[A]$ and let $Q_B=Q[B]$. See Figure \ref{Qpic}.

Note that $d_Q(v)\in\{r, r+1, r+2\}$ for all $v\in V(G)$, since: if $v=x$ it has degree $r+1$ in $Q$; if $v\in\{x_m, y_m\}$ and $n$ is even then it has degree $r$; if $v\in\{y_{m-1}, x_m, y_m\}$ and $n$ is odd then it has degree $r+1$; if $v\in V(M^*)$ and $n$ is odd then it has degree $r+2$, and; otherwise $v$ has degree $r+1$.
Note that for any vertex $v\in V(Q)$, $d_Q(v)=d_A(v)+d_B(v)$, so
by property (P3) this implies that
\begin{equation*}
d_A(v) =\tfrac{1}{2}\left( (d_A(v)+d_B(v)) +(d_A(v)-d_B(v)) \right)\\
\geq \tfrac{1}{2}(r-m^{2/3}).
\end{equation*}
Analogous computations tell us that
\begin{eqnarray}
\delta(Q_A), \delta(Q_B)  &\ge&  \tfrac{1}{2} \left(r- m^{2/3} \right), \label{eqn:minimum-degree-QA-QB} \\
\Delta(Q_A), \Delta(Q_B)  & \le &  \tfrac{1}{2} \left((r+2)+ m^{2/3} \right).  \label{eqn:max-degree-QA-QB}
\end{eqnarray}
Note from our previous observations that $d_Q(x_i)=d_{Q}(y_i)$ for each $i\in [1,m]$. So we have $\sum_{v\in A}d_Q(v)=\sum_{v\in B}d_Q(v)$, and hence
\begin{equation}\label{eqn2}
	|E(Q_A)|=|E(Q_B)|.
\end{equation}
The remainder of  our proof consists of five steps, through which we construct our desired edge-coloring. Before we provide the details for each of these steps, let us first give an outline, along with several necessary definitions. In fact, we shall give two outlines: first an executive summary of the whole process, and then a detailed outline of the five steps for edge-coloring.\\

\noindent \textbf{{\sc Executive Summary of Edge-Coloring Procedure:}}
\begin{itemize}
\item Let $k=\left\lceil \tfrac{1}{2} \left(r+2+ m^{2/3} \right) \right\rceil+4$, and let $\ell=\lceil m^{5/6}  \rceil+1$.
\item In Steps 1-3 we will $k$-edge-color some edges in $Q$. In Steps 1-2 this will  primarily be those in $Q_A, Q_B$, but then in Step 3 we will transform these color classes to be as large as possible, i.e., so that each is a perfect matching.
\item After Steps 1-3 we will not change any of the color classes $1, \ldots, k$, but rather we will add new colors to deal with the remaining edges.
\item In Step 4 we will introduce $\ell$ new colors, and use them to color all the remaining uncolored edges in $Q_A, Q_B$. We will also extend these new color classes to perfect matchings of $Q$ by coloring some additional edges of $Q[A, B]$.
\item In Step 5 we color the remaining necessary edges of $Q$, and again we introduce completely new colors to do so, getting us up to a total of $r+2$ colors. For $n$ even, every edge will be necessary and must be colored. When $n$ is odd, we will finally describe how to delete some edges (between $x$ and $B$ and some edges of $M^*$) that we do not need to color. This fifth and last step is where we employ our Lemma \ref{lem:bipartite-matching-extension}.
\end{itemize}

Step 5 also includes some important details beyond what is written above: because of course our goal is not to edge-color $Q$, but to totally color $G$. In the $n$ even case, we will need to extend the coloring to the two edges incident to $x$, which will then imply our desired total coloring (via the argument given above). In the odd case we will need to describe why it is okay that we didn't color every edge of $Q$, and that it still leads to a total coloring. As the reader will soon see however, it is the edge-coloring of $Q$ (or most of $Q$) that takes nearly all of the effort in this proof.\\

\noindent \textbf{{\sc Detailed Outline of Edge-Coloring Procedure:}}
\begin{itemize}
\item In Step 1 we will give a $k$-edge-coloring $\varphi_0$ to all edges of $Q_A\cup Q_B$ plus a special set $M_1$ of $k$ edges; this whole colored subgraph will be called $Q_{A, B}$. We will ensure that the edges of $M_1$ all get different colors.
\item When $n$ is even, we'll choose $M_1=\{x_1y_1, \ldots, x_ky_k\}$ in Step 1.
\item When $n$ is odd, we'll choose $M_1=\{x_1y_1, \ldots, x_{k-1}y_{k-1}, xx_{m}\}$ in Step 1.
\item If $n$ is even we'll let $\varphi_0=\varphi$ and go directly from Step 1 to Step 3.
\item If $n$ is odd we'll use Step 2 to transform $\varphi_0$ into $\varphi$. The first thing we'll do in Step 2 is say that we may assume $Q, M^*$ have a special property with respect to $\varphi_0$ (involving the colors of $\varphi_0(x)$), subject to uncoloring a small number (at most $m^{2/3}+5$) edges of $Q_A$.
\item The second thing we'll do in Step 2 ($n$ odd only) is to modify the $M_1$ we chose in Step 1, and to change $\varphi_0$ to $\varphi$. In particular, we'll choose four special subsets $M_x^1, M_x^2, M_x^3, M_x^4$ of edges to remove from our initial $M_1$, and instead add to $M_1$ an equal number of edges incident to $x$. The removed edges get uncolored (there will be less than $6m^{2/3}+2\ell$ of them), and the new $M_1$-edges corresponding to the removed edges from $M_x^1$ get colored (the other new $M_1$-edges will already be colored). The new $M_1$ is a subset of $M\cup \delta_Q(x)$, and we ensure it still has exactly $k$ edges, each of which are colored a different color.
\item In Step 3 we will transform $\varphi$ so that each color class is a perfect matching of $Q$. We will do this by pairing up vertices missing a common color $i$ under $\varphi$, and then finding a path of edges between them which alternates between uncolored and $i$-colored edges. We switch on each such path to increase the size of our color class, as desired. We will ensure however that the edges of $M\cup \{xy_{m-1}, xy_m\}$ are not involved in these paths. In particular, this means we'll be maintaining our property that all $k$ edges of $M_1$ are colored a different color (and the other edges of $M$ are uncolored). In this step, some of the edges of $Q_A$ and $Q_B$ will be uncolored.
\item In Step 4 we leave the colors $1, 2, \ldots, k$ alone, but extend $\varphi$ by introducing $\ell$ new colors, say $k+1, \ldots, k+\ell$. We use these new colors to quickly color the rest of the edges in $Q_A, Q_B$ (by Vizing's Theorem), and then are able to extend these further to some edges of $Q[A, B]$ so that each color class is a perfect matching.
\item In Step 5, we define $R$ to be the remaining uncolored edges of $Q$ that we will find necessary to color. When $n$ is even, this is every remaining uncolored edge of $Q$. When $n$ is odd, we will omit from $R$ all the edges in $M_x^1 \cup M_x^3$ (selected in Step 2), as well as any uncolored edges between $X$ and $V(M^*)\cap B$. We then apply Lemma \ref{lem:bipartite-matching-extension} to $R$ (or in the $n$ even case, to $R$ with $x$ and its two incident edges added back) to get our desired edge-coloring. For both $n$ even and $n$ odd we will explain how this implies our total coloring of $G$.
\end{itemize}

Let us now fill in all the details of the above outline.\\

\noindent \underline{\textbf{Step 1:} \emph{Getting a $k$-edge-coloring $\varphi_0$ of $Q_{A,B}$.}}\\

We define $Q_{A,B}=(Q_A \cup Q_B)+M_1$, where $M_1$ is defined as in the above outline (depending on the partity of $n$). Note that in the case where $n$ is odd,
$$|M^*|=\tfrac{1}{2}(r-2) <k-1 = |M_1\cap M|, $$
so by renaming vertices if necessary, we may assume that $M^*  \subseteq M_1\cap M$.

We now apply Lemma \ref{bipartite-matching-extension} to get a $k$-edge-coloring of $Q_{A, B}$ where all the edges of $M_1$  are colored differently. To see that this is possible, first note that $Q_{A, B}$ is a hypergraph of rank $2$, and
$$\Delta(Q_{A, B})+2(2)-1\leq \max\{\Delta(Q_A), \Delta(Q_B)\}+1+3\leq k,$$
where in the last inequality above we used  (\ref{eqn:max-degree-QA-QB}). For the application of Lemma \ref{bipartite-matching-extension}, we use the matching $M_1$, but for the vertex ``$x$'' in the lemma, we simply add a dummy isolate vertex to the graph $Q_{A, B}$, which has no effect on anything (other than vacuously meeting the existence condition in the Lemma). Since $|M_1|=k$, the assumptions of Lemma \ref{bipartite-matching-extension} are satisfied, and we get a $k$-edge-coloring edge coloring $\varphi_0$ of $Q_{A, B}$ where all the edges of $M_1$  receive different colors.

Note that by our careful definition of $Q$, the quantity $|V(Q_{A, B})|$ is always even. So coupling Lemma~\ref{lem:parity} and Lemma~\ref{lem:equitable-coloring-precolored-edges}, we may choose our $k$-edge-coloring $\varphi_0$  such that
\begin{equation}\label{leq4}
\left||\pbar_0^{-1}(i)| -|\pbar_0^{-1}(j)| \right| \le 4
\end{equation}
for any  $i,j\in \{1,\ldots, k\}$.  By renaming some colors if necessary, we may also assume that $\varphi_0(x_iy_i)=i$ for all $i\in \{1,\ldots, k\}$.

For any vertex $v\in V(Q_{A, B})$, $|\pbar_0(v)| =k-d_{Q_{A, B}}(v)$. By substituting into this the value of $k$ and the bound of $\delta( Q_{A, B})\leq \delta(Q_A)$ obtained from (\ref{eqn:minimum-degree-QA-QB}), we get that
\begin{equation}\label{pbar_0v}
|\pbar_0(v)|\leq m^{2/3}+5,
\end{equation}
for every $v\in V(Q_{A, B})$. Since $Q_{A, B}$ has exactly $2m$ vertices, this means that the average number of vertices not seeing any particular colour is at most
$$
\frac{2m(m^{2/3}+5)}{k}  < 4m^{2/3}-4.
$$
By (\ref{leq4}), we get
\begin{equation}\label{eqn1}
	|\pbar_0^{-1}(i)|< 4m^{2/3} \quad \text{for each $i\in \{1,\ldots, k\}$}.
\end{equation}

If $n$ is even we set $\varphi_0=\varphi$ and proceed to Step 3. If $n$ is odd, we transform $\varphi_0$ to $\varphi$ according to Step 2.\\

\noindent \underline{\textbf{Step 2:} \emph{Transforming $\varphi_0$ to $\varphi$ (and modifying $M_1$) when $n$ is odd}.}\\

We would like to say that for every color in $\pbar_0(x)$, that color is used on an edge of $M^*$. Since $|M_1|=k$ and since every edge of $M_1$ is colored differently by $\varphi_0$, every color appears on an edge of $M_1$. However, suppose that some color $i\in \pbar_0(x)$ is used on an edge in $M_1\setminus M^*$. Since $|\pbar_0(x)|  \leq m^{2/3}+5$ (by (\ref{pbar_0v})), and since $|M^*| = \frac{1}{2}(r-2)$, there is some color $j\in \varphi_0(x)$ such that $\varphi(x_jy_j)=j$ and $x_jy_j\in M^*$. But then we can modify $Q, M^*, \varphi_0$ as follows:
\begin{enumerate}[(F1)]
\item delete the edges $xx_j$ and $xy_j$ from $Q$;
\item add the edges $xx_i$ and $xy_i$ to $Q$;
\item delete $x_jy_j$ from $M^*$, and;
\item add $x_iy_i$ to $M^*$.
\item color the edge $xx_i$ with the color $\varphi(xx_j)$, and;
\item if there exists an edge of $Q_A$ that is incident with $x_i$ and is colored by $\varphi_0(xx_j)$, uncolor it.
\end{enumerate}
Note that the above process will uncolor at most $m^{2/3}+5$ edges of $Q_A$. We are okay with this sacrifice, so we may assume we have $Q, M^*, \varphi_0$ where
every color from $\pbar_0(x)$ appears on an edge of $M^*$.

The second half of this step concerns, as described in our outline, modifing the $M_1$ we chose in Step 1, and changing $\varphi_0$ to a final coloring we'll call $\varphi$. We proceed as follows:
\begin{enumerate}[(R1)]
\item let $M_x^1\subseteq M^*$ be the set of edges on which the colors of $\pbar_0(x)$ appear;
\item uncolor all the edges in $M_x^1$, and;
\item color $xy_i$ by the color $i$.
\end{enumerate}
Note that after the above process, all $k$ colors are present at $x$. We continue as follows:
\begin{enumerate}
\item[(R4)] let $M_x^2\subseteq M_1$  be the set of edges on which the colors of $\{\varphi_0(xx_j): j\in \pbar(x)\}$ appear, and;
\item[(R5)] uncolor all the edges in $M_x^2$.
\end{enumerate}
Note that $M_x^1 \cap M_x^2=\emptyset$ since the former involves colors missing at $x$ and the latter involves colors present at $x$.

The sets $M_x^3$ and $M_x^4$ are in a way analogous to $M_x^1, M_x^2$. For the first of these however, only the size really matters. To this end, recall our definition of $\ell= \lceil m^{5/6}  \rceil+1$. We then proceed as follows:
\begin{enumerate}
\item[(R6)] let $M_x^3$ be a set of  $\ell$ edges from $ M^*\setminus (M_x^1\cup M_x^2)$;
\item[(R7)] let $M_x^4$  be the set of  $\ell$ edges of $M_1$ that are colored by colors
from \\$\{\varphi_0(xx_j): \text{$x_jy_j\in M_x^3$ for some $j\in \{1,\ldots, k\}$}\}$, and ;
\item[(R8)] uncolor all the edges in $M_x^3\cup M_x^4$.
\end{enumerate}
In order to see that step (6) is possible, note first that $|M^*\setminus (M_x^1\cup M_x^2)| \geq \tfrac{1}{2}(r-2)-2(m^{2/3}+5)$, and  then consider our assumption about the relative sizes of $r, n$, and the fact that $n$ is sufficiently large.

We claim that we can choose $M_x^3$ and $M_x^4$ so that they are disjoint.
Let $e_1=x_jy_j$ be an arbitrary edge from $M^*\setminus (M_x^1\cup M_x^2)$ for some $j\in \{1,\ldots, k\}$.
Then clearly  $\varphi_0(xx_j) \ne \varphi_0(e_1)$.  We  include $e_1$ in $M_x^3$, and  include the edge of $M_1$ colored by $\varphi_0(xx_j)$
in $M_x^4$.
Suppose for some $t\in \{1,\ldots, \ell\}$, we have chosen $t$ edges for $M_x^3$ and $M_x^4$, respectively, such that   $M_x^3\cap M_x^4=\emptyset$.
Then there are at most $t$ non-chosen edges  from $ M^*\setminus (M_x^1\cup M_x^2)$ that are colored
by colors from $\{\varphi_0(xx_j): \text{$x_jy_j$ is chosen for  $M_x^3$}\}$. Also,
there are at most $t$ unused edges from $x$ to $V( M^*\setminus (M_x^1\cup M_x^2))\cap A$ in $Q$ that are colored
by colors from $\{\varphi_0(x_jy_j): \text{$x_jy_j$ is chosen for $ M_x^3$}\}$.
Then
$$|M^*\setminus (M_x^1\cup M_x^2)|-2t -2t> \tfrac{1}{2}(r-2)-2(m^{2/3}+5)-4\ell \ge 1,$$
where again we use the assumption about the relative sizes of $r, n$, and the fact that $n$ is sufficiently large.
So we can find a non-chosen
edge $e_{t+1}=x_iy_i$,  for some $i\in \{1,\ldots, k\}$,  from $M^*\setminus (M_x^1\cup M_x^2)$
such that  $\varphi_0(x_iy_i)$ is distinct from all the colors on the edges joining $x$
and the vertices from $A$ that are incident with an edge already chosen for $M_x^3$, and that $\varphi_0(xx_i)$
is distinct from all the colors on the chosen edges of $M_x^3$. Now we can add $e_{t+1}$ to $M_x^3$
and add the edge of $M_1$ that is colored by $\varphi_0(xx_i)$ to $M_x^4$.  By the choice of
the two new edges, we know that the updated sets $M_x^3$ and $M_x^4$ are disjoint.
Thus, repeating the above process, we can choose $M_x^3$ and $M_x^4$ such that  $M_x^3\cap M_x^4=\emptyset$.

Finally, with all of the above done, update $M_1$ to be  $M_1\setminus (M_x^1\cup M_x^2 \cup M_x^3 \cup M_x^4)$ union
the set of edges that are incident with $x$ and are colored by a color that was initially assigned to an edge from $M_x^1\cup M_x^2 \cup M_x^3 \cup M_x^4$. This adjustment means that the set $M_1$ still contains exactly one edge of each of the $k$ colors, but now these include more edges of $\delta_Q(x)$.
Note  that the overall effect of the above process was to color more of the edges in $\delta_Q(x)$, and uncolor at most $|M_1^x \cup M_2^x \cup M_x^3 \cup M_x^4| \le  2 (m^{2/3}+5)+2\ell$ edges. \\

\noindent \underline{\textbf{Step 3:} \emph{Transforming the $k$ color classes from $\varphi$ into perfect matchings.}}

During this step, by swapping along alternating paths, we will increase the size of the $k$ color classes obtained in Step 1 until each color class is a perfect matchings of $Q$. During this process, we will maintain
the colors on the edges of $M_1$ and we will maintiain the colors on the edges incident with $x$ (when $n$ is odd). We will however, uncolor some of the edges of $Q_A$ and $Q_B$.
%nd will color some of the edges of $H-E(Q_{A,B})-M_2$.
Let $R_A, R_B$ be the the subgraphs induced by all uncolored edges in $Q_A, Q_B$, respectively. These graphs will both initially be empty if $n$ is even, and $Q_B$ will also be initially empty even if $n$ is odd. As we progress through this step however,
one or two edges will be added to each of $R_A$ and $R_B$
each time we swap colors on an alternating path. We will tighlty control these uncolored edges, and will ensure that the following conditions are satisfied after the completion of Step 3:
\begin{enumerate}
	\item[(C1)] $|E(R_A)|=|E(R_B)| <4m^{5/3}$.
		\item[(C2)] $\Delta(R_A),\Delta(R_B) < m^{5/6}$.
	\item[(C3)] Each vertex of $Q$ is incident with fewer than $2m^{5/6}$ colored edges of $Q[A, B]$, i.e., the bipartite subgraph of $Q$ induced by the bipartition $(A, B)$.
\end{enumerate}

At any point during our process, we say that an edge  $e=uv$ is \emph{good} if $e\not \in E(R_A) \cup E(R_B)$,  the degree of $u$
and $v$ in both $R_A$ and  $R_B$ is less than $m^{5/6}-1$, and in the case that $n$ is odd, $e$ is not incident with any vertex from $V(M_x^3)\cap B$. In particular, a good edge can be added to $R_A$ or $R_B$ without violating condition (C2).

By Lemma \ref{lem:parity}, we know that $|\pbar^{-1}(i)|$ is even for every color $i\in \{1,\ldots, k\}$, so since $(A, B)$ is a partition of $V(Q)$, the quantity $|\pbar^{-1}_A(i)|-|\pbar^{-1}_B(i)|$ is also even. We can pair up as many vertices as possible from $\pbar^{-1}_A(i),\pbar^{-1}_B(i)$, and then pair up
the remaining unpaired vertices from  $\pbar^{-1}_A(i)$ or $\pbar^{-1}_B(i)$.
We call each of these pairs a \emph{missing-common-color pair} or \emph{MCC-pair} with respect to the color $i$. In fact, given any MCC-pair $(a,b)$ with respect to some colour $i$, we will show that the vertices $a, b$ are joined by some path $P$ (with at most 7 edges), where one end of the path is an uncolored edge from $Q[A, B]$, and then the path alternates between good edges (colored $i$) and uncolored edges. Given such a path $P$, we will swap along $P$ by giving all the uncolored edges $i$, and uncoloring all the colored edges. After this, both $a$ and $b$ will be incident with edges colored $i$, and at most three good edges will be added to $R_A\cup R_B$. Given the alternating nature of the path $P$, at most two good edges will be added to any one of $R_A, R_B$.
The main difficulty in this Step 3 (and perhaps in this proof) is to show that such path $P$ exists. Before we do this, let us show that conditions (C1)-(C3) can be guaranteed at the end of this swapping process, that is, at the end of Step 3. In fact,
we will only ever add good edges to and $R_A$ and $R_B$, so condition (C2) will  hold automatically.  It remains then only to check (C1) and (C3).

Let us first consider condition (C1).  We deleted at most $m^{2/3}+5$ edges of $Q_A$ that were incident with $x$ after Step 1 but we have $\pbar_0(x)=\pbar(x)$, and we
uncolored  at most $m^{2/3}+5$ edges of  $Q_A$ but for
each of these edges, one of its endvertices has the same set of missing colors under both $\varphi_0$ and $\varphi$.  So the sum of
missing colors  under $\varphi$ taken over all vertices of $Q$ is at most $2m(m^{2/3}+5)+2(m^{2/3}+5)<3m^{5/3}$.
Thus there are less than  $1.5m^{5/3}<2m^{5/3}$ MCC-pairs.

For each MCC-pair $(a,b)$ with $a,b\in V(Q)$,
at most  two edges will be added to each of $R_A$ and $R_B$
when we exchange an alternating path from $a$ to $b$.  Thus
there will always be less than
$4m^{5/3}$  uncolored
edges in each of $R_A$ and $R_B$.
At the completion of
Step 3, each of the $k$ color classes is a perfect matching of $Q$
so each of $Q_A$ and $Q_B$ have the same number of
colored edges. Since $|E(Q_A)|=|E(Q_B)|$ by~\eqref{eqn2},
we have $|E(R_A)|=|E(R_B)|$. As $R_B$ is initially empty and less than $4m^{5/3}$ edges will be added to $R_B$
by the  end of Step 3, we have $e(R_B) <4m^{5/3}$. As $|E(R_A)|=|E(R_B)|$,
 we also have    $e(R_A) <4m^{5/3}$. Thus  Condition (C1) will be satisfied at the end of Step 3.

We  now show that Condition (C3)  will also be satisfied at the end of Step 3. In the process of Step 3,   the number of newly colored edges of $Q[A, B]$ that are incident with a vertex   $u\in V(Q)$ will equal  the number of our chosen alternating paths containing $u$. The number of such alternating paths of which $u$ is the first vertex will equal the number of colors missing from $u$ at the end of Step 2.  The number of alternating paths in which $u$ is not the first  vertex will  equal  the degree of $u$ in $R_A\cup R_B$, and so will be less than $m^{5/6}$.   As the edges in $M_1$ were
colored in Steps 1-2,
the number of colored edges of $Q[A, B]$ that are incident with $u$ will be less than $|\pbar(u)|+m^{5/6}+1<2m^{5/6}$, where note that   $|\pbar(u)|=k-|\varphi(u)| <3 m^{2/3}$  as we uncolored at most $k-\frac{1}{2}(r-m^{2/3})<m^{2/3}+4$ edges in $Q_A$
when we dealt with the colors in $\pbar_0(x)$ when $n=2m-1$,
and all those edges could be all incident with one single vertex in $Q_A$.
Thus Condition (C3)   will be satisfied.

Let $i\in\{1,  \ldots, k\}$. We will describe how to find the alternating paths for the MCC-pairs with respect to $i$. There is one $i$-colored edge in $M_1$, say $e_i$ with endpoints $V(e_i)$, and we will be sure to avoid it in all of the alternating paths joining our MCC-pairs. For each vertex in an MCC-pair, we will see how to construct a short (two-edge) alternating path from it that avoids $e_i$. We will see that there are so many choices for such paths (even avoiding the vertex set $M_x^3$ as we will wish to do), that we can connect them to link each MCC-pair.

Let $v\in V(Q)$.
Suppose that $v\in S_v\in\{A, B\}$, and let $T_v$ be the other one of $A, B$. Define $N_1(v)$ to be the set of all vertices in $T_v\setminus (V(e_i)\cup  V(M_x^3) \cup\{x\})$, where $V(e_i)$ is the set of the two endvertices of $e_i$,  that are joined to $v$ by an uncolored edge, and are incident to a good edge colored $i$.  Let $N_2(v)$ be the set of vertices in $T_v\setminus (V(e_i)\cup  V(M_x^3) \cup\{x\})$ that are joined to a vertex of $N_1(v)$ by a good edge of color $i$. Note that we have $|N_1(v)|=|N_2(v)|$, but some vertices may be in $N_1(v)\cap N_2(v)$.

There are fewer than $4m^{5/3} $ edges in $R_B$ (by (C1)), so there are fewer than $8m^{5/6}$ vertices of degree at least $m^{5/6}$ in $R_B$.
Each non-good edge is incident with one or two
vertices of $R_B$ through the color $i$. Thus   the number of vertices in $B$ incident with a non-good edge colored by $i$
is less than  $16m^{5/6}$.
In addition, by~\eqref{eqn1} there are fewer than $3m^{2/3}$
vertices in $B$ that are missed by the color $i$.  There are
at most $2|V(M_x^3)\cap B|$ vertices of $B$ that are incident with an edge colored by $i$ such that one endvertex
of the edge is from $V(M_x^3)\cap B$.
So the number of vertices in $B\setminus\{y_i\}$ that are
incident with a non-good edge colored $i$ is less than
\begin{equation}\label{eqn3}
	16m^{5/6}+3m^{2/3}+2 \ell<19m^{5/6}.
\end{equation}
By symmetry,
the number of vertices  in $A\setminus\{x_i\} $ that are
incident with a non-good edge colored $i$ is less than
$19m^{5/6}$. Therefore for any $v\in V(Q)$, we have
\begin{eqnarray}
|N_1(v)|=|N_2(v)|  & \ge& d_{Q[A, B]}(v) -2-2m^{5/6} -19 m^{5/6}-2  \nonumber  \\
 &\ge & \tfrac{1}{2} \left((1+\ve)n- m^{2/3} \right)-22m^{5/6} \nonumber \\
 &>&\tfrac{1}{2}(1+0.5\ve) m, \label{eqn:N-M-sizea}.
\end{eqnarray}
where  in the inequalities to get~\eqref{eqn:N-M-sizea}, the first term $-2$ is to make sure that the endpoints of $e_i$ are not counted,  the number $2m^{5/6}$ is the maximum number of
colored edges of $Q[A, B]$ that are incident with $v$ at the end of Step 3 (as per (C3)), and the last term $-2$ is to make sure that when $n$ is odd, the vertex $x$ and another vertex that is adjacent with $x$ through a good edge colored by $i$ is not counted.

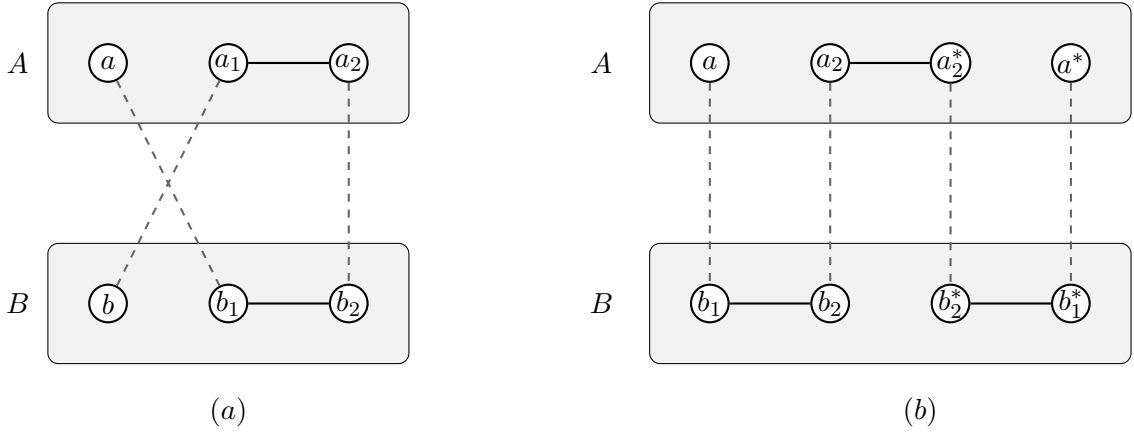
\begin{figure}[!htb]
	\begin{center}
		
		\begin{tikzpicture}[scale=0.8]

			\begin{scope}[shift={(0,0)}]
				\draw[rounded corners, fill=white!90!gray] (0, 0) rectangle (6, 2) {};
				
				\draw[rounded corners, fill=white!90!gray] (0, -4) rectangle (6, -2) {};
				
				{\tikzstyle{every node}=[draw ,circle,fill=white, minimum size=0.5cm,
					inner sep=0pt]
					\draw[black,thick](1,1) node (a)  {$a$};
					\draw[black,thick](3,1) node (a1)  {$a_1$};
					\draw[black,thick](5,1) node (a2)  {$a_2$};
					\draw[black,thick](1,-3) node (b)  {$b$};
					\draw[black,thick](3,-3) node (b1)  {$b_1$};
					\draw[black,thick](5,-3) node (b2)  {$b_2$};
				}

				\path[draw,thick,black!60!white,dashed]
				(a) edge node[name=la,pos=0.7, above] {\color{blue} } (b1)
				(a2) edge node[name=la,pos=0.7, above] {\color{blue} } (b2)
				(b) edge node[name=la,pos=0.6,above] {\color{blue}  } (a1)
				;
				
				\path[draw,thick,black]
				(a1) edge node[name=la,pos=0.7, above] {\color{blue} } (a2)
				(b1) edge node[name=la,pos=0.7, above] {\color{blue} } (b2)
				;
				
				\node at (-0.5,1) {$A$};
				\node at (-0.5,-3) {$B$};
				\node at (3,-4.8) {$(a)$};

			\end{scope}

			\begin{scope}[shift={(2,0)}]
				\draw[rounded corners, fill=white!90!gray] (8, 0) rectangle (16, 2) {};
				
				\draw[rounded corners, fill=white!90!gray] (8, -4) rectangle (16, -2) {};
				
				{\tikzstyle{every node}=[draw ,circle,fill=white, minimum size=0.5cm,
					inner sep=0pt]
					\draw[black,thick](9,1) node (c)  {$a$};
					\draw[black,thick](11,1) node (c1)  {$a_{2}$};
					\draw[black,thick](13,1) node (c2)  {$a^*_{2}$};
					\draw[black,thick](15,1) node (c3)  {$a^*$};
					%\draw[black,thick](12,1) node (c2)  {$a_{12}$};
					%\draw[black,thick](13.5,1) node (c3)  {$a_3$};
					%\draw[black,thick](15,1) node (c4)  {$a_2$};
					\draw[black,thick](9,-3) node (d)  {$b_{1}$};
					\draw[black,thick](11,-3) node (d1)  {$b_{2}$};
					\draw[black,thick](13,-3) node (d2)  {$b^*_{2}$};
					\draw[black,thick](15,-3) node (d3)  {$b^*_{1}$};
					%\draw[black,thick](12,-3) node (d2)  {$b_2$};
					%\draw[black,thick](13.5,-3) node (d3)  {$b_{22}$};
					%\draw[black,thick](15,-3) node (d4)  {$b_{21}$};
					
				}
				\path[draw,thick,black!60!white,dashed]
				(c) edge node[name=la,pos=0.7, above] {\color{blue} } (d)
				(c1) edge node[name=la,pos=0.7, above] {\color{blue} } (d1)	
				(c3) edge node[name=la,pos=0.7, above] {\color{blue} } (d3)	
				(c2) edge node[name=la,pos=0.7, above] {\color{blue} } (d2)	
				;
				
				\path[draw,thick,black]
				(d) edge node[name=la,pos=0.7, above] {\color{blue} } (d1)
				(c2) edge node[name=la,pos=0.7, above] {\color{blue} } (c1)
				(d3) edge node[name=la,pos=0.7, above] {\color{blue} } (d2)
				;
				\node at (7.2,1) {$A$};
				\node at (7.2,-3) {$B$};
				\node at (12.5,-4.8) {$(b)$};	
			\end{scope}	
			
		\end{tikzpicture}
	\end{center}
	\caption{The alternating path $P$. Dashed lines indicate uncoloured edges, and solid
		lines indicate edges with color $i$.}
	\label{f1}
\end{figure}

Suppose now that we have some MCC-pair $(a,b)$ with respect to $i$, with $a\in A$ and $b\in B$.
By~\eqref{eqn:N-M-sizea}, we have
$
|N_2(a)|,  |N_2(b)|>\tfrac{1}{2}(1+0.5\ve) m.
$
We choose $a_1a_2$ with color $i$ such that $a_1\in N_1(b)$
and $a_2\in N_2(b)$. Now as $|N_1(a)|, |N_1(a_2)| >\tfrac{1}{2}(1+0.5\ve) m$ by~\eqref{eqn:N-M-sizea}, we know that $|N_1(a_2)\cap N_2(a) |>\frac{1}{2}\ve m$. We choose $b_2\in N_1(a_2)\cap N_2(a)$
so that  $b_2a_2 \not\in (M\cup \{xy_{m-1}, xy_m\})\setminus M_1$,
and let $b_1\in N_1(a)$ such that $b_1b_2$ is colored by $i$.
Then $P=ab_1b_2a_2a_1 b$ is an alternating path from $a$ to $b$ (See Figure~\ref{f1}(a)). We exchange $P$ by coloring $ab_1, b_2a_2$ and $a_1b$
with color $i$ and uncoloring the edges $a_1a_2$ and $b_1b_2$.
After the exchange,  the color $i$ appears on edges incident with $a$ and $b$,
the edge $a_1a_2$ is added to $R_A$
and the edge $b_1b_2$ is added to $R_B$.

Consider then an MCC-pair $(a,a^*)$
with respect to  $i$ such that $a,a^* \in A$ (the case for an MCC-pair $(b,b^*)$
with $b,b^* \in B$ can be handled similarly).
By~\eqref{eqn:N-M-sizea}, we have  $|N_2(a^*)|>\tfrac{1}{2}(1+0.5\ve) m$.
We take an edge $b_1^*b_2^*$ colored by $i$ with $b_1^* \in N_1(a^*)$ and  $b_2^*\in N_2(a^*)$.
Then again, by~\eqref{eqn:N-M-sizea}, we have
$
|N_2(a)|,  |N_2(b_2^*)|>\tfrac{1}{2}(1+0.5\ve) m.
$
Therefore, as each vertex  $c\in N_2(b_2^*)$
satisfies $|N_1(c)|>\tfrac{1}{2}(1+0.5\ve) m$,  we have $|N_1(c)\cap N_2(a)| >\tfrac{1}{2}\ve m$. We take $a_2a_2^*$ colored by $i$ with $ a_2^*\in N_1(b_2^*)$ such that $a_2^*b_2^*\not\in (M\cup \{xy_{m-1}, xy_m\})\setminus M_1$ and  $a_2\in   N_2(b_2^*)$.
Then we let $b_2\in N_1(a_2)\cap N_2(a)$ such that $a_2b_2\not\in (M\cup \{xy_{m-1}, xy_m\})\setminus M_1$, and let $b_1$
be the vertex in $N_1(a)$ such that $b_1b_2$ is colored by $i$.
Now we get the alternating path $P=ab_1b_2 a_2 a_2^* b_2^* b_1^* a^*$ (See Figure~\ref{f1}(b)).
We exchange $P$ by coloring $ab_1, b_2a_2, a_2^*b_2^*$ and $b_1^*a^*$
with color $i$ and uncoloring the edges $b_1b_2, b_1^*b_2^*$ and $a_2a_2^*$.
After the exchange,  the color $i$ appears on edges incident with $a$ and $a^*$,
the edges $b_1b_2$ and $b_1^*b_2^*$ are added to $R_B$
and the edge $a_2a_2^*$ is added to $R_A$.  In this process,
we  added  one edge  to $R_A$
and at most two edges to $R_B$.

We have shown above how to find alternating paths between all MCC-pairs, and to modify $\varphi$ by switching on these.  Each time we make such a switch, we increase by one the size of the color class $i$ (and decrease by 2 the number of vertices that are missing color $i$). By repeating this process, we can therefore continue until the color class $i$ is a perfect matching of $Q$. Durig this whole process, note that we did not alter any edge-colorings for edges from $M_1$, nor any edges from $M\cup \{xy_{m-1}, xy_m\}$. Hence we still have that all edges in $M_1$ receive different colors and all edges in $(M\cup \{xy_{m-1}, xy_m\})\setminus M_1$ are uncolored. \\

\noindent \underline{\textbf{Step 4:} \emph{Coloring $R_A$ and $R_B$ and extending the new color classes.}}\\

In this step we leave the colors $1, 2, \ldots, k$ alone, but extend $\varphi$ by introducing $\ell$ new colors, say $k+1, \ldots, k+\ell$. By (C2), $\Delta(R_A), \Delta(R_B)<m^{5/6}$. Recalling that $\ell= \lceil m^{5/6}  \rceil+1$, Vizing's Theorem gives us an $\ell$-edge-coloring of $R_A\cup R_B$ using the colors $k+1, \ldots, k+\ell$.  By Theorem~\ref{lem:equa-edge-coloring}, we may assume that all these new color classes differ in size by at most one. Since $|E(R_A)|=|E(R_B)|$ by (C1),
by possibly renaming some color classes we may assume that each color
appears on the same number of edges in $R_A$ as it does in $R_B$.

By (C1) we get that $|E(R_A)|, |E(R_B)|< 4m^{5/3}$, so since  $\ell>m^{5/6}$, each of
the colors $k+1,\ldots, k+\ell$ appears on fewer than $4 m^{5/6} +1$
edges in each of $R_A$ and $R_B$. We will now color some of the edges of $G[A, B]- M$ with the $\ell$ new colors so that each class induces a perfect matching. To this end, we
 perform the following procedure for each of the $\ell$ colors in turn.

Let $i$ be a color in $\{k+1,\ldots, k+\ell\}$. We define $A_i, B_i$ as the sets of vertices in $A, B$, respectively, that are incident with edges colored $i$.
Then, from our discussion in the previous paragraph, $|A_i|,|B_i| < 2 (4 m^{5/6} +1)$.   By our choices we in fact have $|A_i|=|B_i|$;
let $H_i$ be the subgraph of $Q[A, B]-M$ obtained by
deleting the vertex sets $A_i \cup B_i$ and removing all colored edges. We will show that $H_i$ has a perfect matching and we will color these matching edges with $i$.

By (C3), each vertex in $V(Q)$ is incident with fewer than
$2m^{5/6}+\ell \le 4 m^{5/6}$
edges of $Q[A, B]$ that are colored. Also each vertex in $A$ has fewer than $ 2(4 m^{5/6} +1)$ edges to $B_i$ and each vertex in $B$ has fewer than $ 2(4 m^{5/6} +1)$ edges to $A_i$. So each vertex  $v\in V(H_i)$ (with $v\ne x$ when $n$ is odd) is adjacent in $H_i$ to  more than
$$
\tfrac{1}{2}\left((1+\ve)m-m^{2/3}\right)-4m^{5/6}-2(4 m^{5/6} +1)-2>\tfrac{1}{2}(1+0.5\ve)m
$$
vertices, where the term $-2$ is due to the one or two edges in $M$ that are possibly incident with $v$ in $Q$.

We are now ready to extend our color classes. If $n$ is even, then $H_i$
has a perfect matching $F$ by Lemma~\ref{lem:matching-in-bipartite};
%since $\delta(H_i)\geq |V(H_i)|/2$;  No the degree is about $|V(H_i)|/4$.
we color the edges of $F$ with the color $i$ to get our desired result. So we may assume that $n$ is odd.
Recall from Step 2 that $|M_x^3|=\ell$ and we have $V(M_x^3)\cap V(R_B) =\emptyset$.
Also  $x$ is not contained in $R_A$
by our construction (see~\eqref{eqn:N-M-sizea}).   Thus all the  edges in $E_{Q}(x, (V(M^*)\setminus V(M_x^1)) \cap B)$ are uncolored in Step 2,
and that  $\{x\} \cup (V(M_x^3)\cap B) \subseteq  V(H_i)$. Thus, for each $i\in \{k+1,\ldots, k+\ell\}$,
we choose a distinct edge, say $xy$ from $E_Q(x, V(M_x^3)\cap B)$ and color it by $i$.
Then $H_i-\{x,y\}$
has a perfect matching $F$ by virtue of its minimum degree (as in the even case), and we get our desired extension. \\

\noindent \underline{\textbf{Step 5:} \emph{ Coloring the remaining graph $R$ and getting a total coloring of $G$.}}\\

As described in the outline, we define $R$ to be the remaining uncolored edges of $Q$ that we will find necessary to color. When $n$ is even, this is every remaining uncolored edge of $Q$. When $n$ is odd, we will omit from $R$ all the edges in $M_x^1 \cup M_x^3$ (selected in Step 2), as well as any uncolored edges between $x$ and $V(M^*)\cap B$.

Note that when $n$ is even,  $\Delta(R) =r+1-k-\ell$.
When $n$ is odd, we need to make a few observations before being able to bound the maximum degree. First, only vertices in $V(M^*)\cap B$  have degree $r+2$ in $G^*$, with all others having degree $r+1$ (since $n$ is odd). So it follows that for any $w\in B $, we have $r-k-\ell  \le d_{R}(w)  \le r+1-k-\ell$;
and for any $w\in A$, we have $r-k-\ell  \le d_{R}(w)  \le r+2-k-\ell$.
In any case, we have
$$
 \Delta(R)     \ge  r-k-\ell
 \ge  r-\left(\left\lceil\tfrac{1}{2} \left(r+ m^{2/3} \right) \right\rceil+4\right) -( \lceil m^{5/6}  \rceil+1)
 >  \tfrac{1}{2} (1+0.5\ve) m. $$

In the case where $n$ is even, we let  $R^*=R+\{xx_{m}, xy_{m}\}$, that is the remaining edges of $Q$ plus the last two edges of $G^*$ that are needed to achieve our total coloring.  If $n$ is odd, we just let $R^*=R$. We now show that we can apply Lemma~\ref{lem:bipartite-matching-extension} to $R^*$.

Since all edges of $Q_A, Q_B$ were colored in previous steps, $R^*$ is indeed a bipartite graph, with a matching $M_2:=M\setminus M_1$), plus an additional vertex $x$ with edges from $x$ to $V(G)\setminus V(M)$ when $n$ is even, and from $x$ to $\{y_{m-1}, y_m\}$ when $n$ is odd.  As $|M|  \le  m-1$ and $M_1$ and edges  in $E_Q(x, V(M_x^1 \cup M_x^3))$ covers  $2k-2|M_x^2|-2|M_x^4|$ vertices of $V(M)$, we get
$$|M_2| \le m-1-k+(m^{2/3}+5)+\ell<\frac{1}{2}m.$$

When $n$ is even, we have $\Delta(R^*)=r+1-k-\ell$.
As $d_{R^*}(x)=2$,
we have $ \Delta(R^*)+1 \ge |M_2| +|\delta_{R^*}(x)|$, and so
the assumptions of Lemma~\ref{lem:bipartite-matching-extension}
is satisfied. So $R^*$
has a $(\Delta(R^*)+1)$-edge-coloring where all the edges of $M_2\cup \delta_{R^*}(x)$ are colored differently.

When $n$ is odd, we have $\Delta(R^*)=r+2-k-\ell$.
As all the edges in $M_x^1 \cup M_x^2 \cup E_Q(x,  \big(V(M^*)\setminus (V(M_x^1 \cup M_x^2))\big)\cap B)$ were removed when we get $R^*$ from $R$, it follows that $d_{R^*}(w) \le r+1-k-\ell$ for every $w\in B$.
As $d_{R^*}(x)=2$,
we have $ \Delta(R^*) \ge |M_2| +|\delta_{R^*}(x)|$, and so
the assumptions of Lemma~\ref{lem:bipartite-matching-extension}
is satisfied. So $R^*$
has a $\Delta(R^*)$-edge-coloring where all the edges of $M_2\cup \delta_{R^*}(x)$ are colored differently.

When $n$ is even, combining the color classes obtained above in $R^*$ with the perfect matchings obtained in Steps 3 and 4, we have obtained an edge coloring  $\varphi^*$ of $G^*$ using
$r+2$ colors such that all the edges in $M\cup \delta_{G^*}(x)$ are colored differently. So, we get our desired total coloring of $G$ using $(r+2)$ colors, as described early in our proof.

Now suppose that $n$ is odd. Combining the color classes in $R^*$ with the perfect matchings obtained in Steps 3 and 4, we have obtained an edge coloring of
$$G^{**}:=G^*-M_x^1 -M_x^2-E_{G^*}(x, (V(M^*)\setminus (V(M_x^1) \cup V(M_x^3))) \cap B).$$
This edge coloring uses $r+2$ colors and all the edges in
$$M^{**}:=(M\setminus (M_x^1 \cup M_x^3)) \cup E_{G^*}(x, \{x_{m-1}, y_{m-1}, y_{m}\} \cup V(M_x^1) \cup V(M_x^3))$$
are colored differently.
Note that $M^{**}$ covers all the vertices of $G$.  Define $\psi:V(G)\cup E(G) \rightarrow \{1, 2, \ldots, r+2\}$ by: $\psi(e)=\varphi^*(e)$ for all $e\in E(G)$;
$\psi(u)=\varphi^*(e)$ for all $e\in M^{**}$ such that $u\in V(G)$ and $u$ is incident with $e$ in $G^*$.   By the same argument as in the even case, $\psi$
is a total coloring of $G$ using at most $r+2$ colors.
This completes  the proof of Theorem~\ref{thm:main0}.
\qed

\bibliographystyle{abbrv}
\bibliography{BIB}
\end{document}